\documentclass[12pt,twoside,final,psamsfonts]{amsart}
\usepackage[draft=false]{hyperref}
\usepackage[psamsfonts]{amssymb}
\usepackage{times,a4wide}

\theoremstyle{plain}
\newtheorem*{theorem*}{Theorem}
\newtheorem*{theoremA}{Theorem A}
\newtheorem{theorem}{Theorem}
\newtheorem{proposition}[theorem]{Proposition}
\newtheorem*{proposition*}{Proposition}

\newtheorem*{corollary*}{Corollary}
\newtheorem{lemma}[theorem]{Lemma}
\newtheorem*{lemma*}{Lemma}
\newtheorem{claim}[theorem]{Claim}
\theoremstyle{definition}
\newtheorem{remark}[theorem]{Remark}
\newtheorem*{remark*}{Remark}
\newtheorem*{conjecture*}{Conjecture}
\theoremstyle{definition}
\newtheorem*{definition}{Definition}

\newcommand{\C}{\mathbb{C}}
\renewcommand{\H}{\mathbb{H}}
\newcommand{\N}{\mathbb{N}}
\newcommand{\R}{\mathbb{R}}
\newcommand{\E}{\mathcal{E}}
\newcommand{\EP}{\hat{\E^\prime}(\R)}

\newcommand{\Har}{\operatorname{Har}}
\newcommand{\lil}{\lambda\in\Lambda}
\renewcommand{\Im}{\operatorname{Im}}
\renewcommand{\Re}{\operatorname{Re}}

\title[Interpolation in $\hat{\E^\prime}(\R)$]
{A Geometric Characterization of Interpolation in $\hat{\E}^\prime(\R)$}

\author[X. Massaneda, J. Ortega-Cerd\`a, M. Ouna\"{\i}es]
{Xavier Massaneda, Joaquim Ortega-Cerd\`a \& Myriam Ouna\"{\i}es}

\address{Departament de Matem\`atica Aplicada i An\`alisi,
Universitat  de Bar\-ce\-lo\-na, Gran Via 585, 08071-Bar\-ce\-lo\-na, Spain}
\email{xavier@mat.ub.es}

\address{Departament de Matem\`atica Aplicada i An\`alisi,
Universitat  de Bar\-ce\-lo\-na, Gran Via 585, 08071-Bar\-ce\-lo\-na, Spain}
\email{jortega@ub.edu}

\address{Institut de Recherche Math\'ematique Avan\c c\'ee, Universit\'e 
Louis Pasteur 7 Rue Ren\'e Des\-car\-tes, 67084 Strasbourg CEDEX, France.}
\email{ounaies@math.u-strasbg.fr}

\thanks{First and second authors supported by DGICYT grant 
BFM2002-04072-C02-01 and the CIRIT grant 2001-SGR00172.}
\date{\today}



\begin{document}

\begin{abstract} 
We give a geometric description of the interpolating varieties for the algebra
of Fourier transforms of distributions (or Beurling ultradistributions) 
with compact support on the real line.
\end{abstract}

\maketitle

\section{Introduction}
Let $\E(\R)$ be the space of smooth functions in $\R$ and let $\E^\prime(\R)$
be its dual, the space of distributions with compact support on $\R$. It is
well known that the space $\hat{\E^\prime}(\R)$ of Fourier transforms of
distributions in $\E^\prime(\R)$ coincides with the algebra  of entire
functions $f$ such that 
\[
|f(z)|\leq C (1+|z|)^A e^{B |\Im\, z|},
\] 
where $A, B, C>0$ may depend on $f$ (see \cite[Theorem 1.4.15]{Br-Ga}).

A discrete sequence $\Lambda\subset \C$ is called $\EP$-interpolating 
when the interpolation problem $f(\lambda)=v_\lambda$, $\lil$, has a solution
$f\in\EP$ for every sequence of complex values $\{v_\lambda\}_{\lil}$
having the characteristic growth of $\EP$ on $\Lambda$ (see the precise
definition in Section 2).

The origin of the interest in $\EP$-interpolation lies in its
relationship with convolution equations and, in particular, with the 
density of exponential families $\{e^{i\lambda x}\}_{\lil}$ in the
space of solutions $g\in\mathcal E(\R)$ of equations of type
$\mu\star g=0$, $\mu\in\mathcal E^\prime(\R)$. Any solution $g$ 
to the convolution equation is the limit of linear combinations of 
$\{e^{i\lambda x}\}_{\lil}$ where $\Lambda$ is the zero set of $\hat \mu$. If moreover
the sequence $\Lambda$ is $\EP$-interpolating then the series that represents
$g$ enjoys better converging properties. For more on this relationship 
see \cite{Eh-Ma} or \cite[Chapter 6]{Br-Ga} (in particular Theorem 6.1.11).

The space $\EP$ is a particular case of the algebras 
\[
A_p=\{f\in H(\C) : \log|f(z)|\leq A +B p(z)\ \textrm{for some}\ A,B>0\}
\]
associated to 
positive measurable weights $p$, 
obtained by taking 
\begin{equation}\label{pes}
p(z)=|\Im\, z|+\log(1+|z|^2).
\end{equation}
There exists an analytic characterization of interpolating sequences for
general $A_p$ spaces when $p$ is subharmonic  (see Theorem A below). However, a
complete geo\-metric description exists only for subharmonic wei\-ghts $p$
which are both radial $(p(z)=p(|z|))$ and doubling ($p(2z)\leq C p(z)$ for some
$C>0$); see
\cite[Corollary 4.8]{Br-Li}.

For the weight \eqref{pes} Ehrenpreis and Malliavin gave a necessary geometric
condition  which turns out to be sufficient  provided that $\Lambda$  is a zero
sequence of a slowly decreasing function  (see  \cite[Theorem 4]{Eh-Ma} and its
proof). Later Squires, probably unaware of Ehrenpreis and Malliavin's result
(which was stated in terms of solutions to convolution equations), proved the
same result \cite[Theorem 2]{Sq}.

In this paper we give a geometric characterization for
$\EP$-inter\-po\-la\-ting  sequences (Theorem~\ref{teorema}). The
characterization shows in particular that the geometric condition given by
Ehrenpreis \& Malliavin and Squires is also sufficient whenever the sequence is
contained in the region  
\[
|\Im\, z|\leq C\log(1+|z|^2).
\] 
In general, however, their condition alone is not sufficient.

A similar characterization is obtained for the more general Beurling weights,
that is, those of the form $p(z)=|\Im\, z|+\omega(|z|)$ where
$\omega:[0,\infty)\longrightarrow [0,\infty)$ is normalized  with
$\omega(0)=0$, subadditive, with $\log(1+t)\lesssim \omega(t)$ for $t>1$  and
$\int_0^\infty \frac{\omega(t)}{1+t^2}\, dt<\infty$.  Beurling weights appear
naturally in the context of convolution equations when one replaces
distributions with compact support with  Beurling-Bj\"orck ultradistributions
of compact support (see \cite{Bjorck66}). These weights  are not necessarily
subharmonic, but we will prove that they are equivalent to a subharmonic weight
(see Lemma~\ref{regularization}).  

The paper is structured as follows. In the following section we give the
precise definition of interpolating variety, introduce the background on the
problem and  state the main result. In Section 3 we prove that the geometric
conditions of  Theorem~\ref{teorema} are necessary, while in Section 4 we show
that they are  also sufficient.

A final remark about notation.  $C$ will always denote a positive constant and 
its actual value may change from one  occurrence to the next. $A= O(B)$ and 
$A\lesssim B$ mean that $A\leq c B$ for some $c>0$, and $A\simeq B$ is
$A\lesssim B\lesssim A$.

\section{Preliminaries}\label{preliminaries}
For the following definition and general background on the problem
we refer to 
\cite[Chapter~2]{Br-Ga}.

A measurable function  $p:\C\longrightarrow\R_+$ is called a \emph{weight}
if for some $C,D>0$:
\begin{itemize}
\item[(a)] $ \log(1+|z|^2)\leq C p(z) $ for all $z\in\C$.
\item[(b)] $p(\zeta)\leq Cp(z)+D\; $ if $\; |\zeta-z|\leq 1$. 
\end{itemize}
The importance of these properties lies in the consequences on the ring $A_p$
defined  in the introduction: (a) implies that $A_p$ contains all polynomials,
and (b) that $A_p$ is closed under differentiation.

The algebra $A_p$ can be thought of as the union of the Hilbert spaces 
\[
A_{p,\alpha}^2:=\Bigl\{f\in H(\C) : \|f\|_{A_{p,\alpha}^2}^2= \int_{\C}
|f(z)|^2 e^{-\alpha p(z)}dm(z)<\infty\Bigr\} 
\]
for $\alpha>0$, as well as the union of the Banach spaces
\[
A_{p,\alpha}^\infty:=\Bigl\{f\in H(\C) : \|f\|_{A_{p,\alpha}^\infty}
=\sup_{z\in\C} |f(z)| e^{-\alpha p(z)}<\infty\Bigr\}\ .
\]
Also, $A_p=\cup_{\alpha>0}A_{p,\alpha}^\infty$ has structure of (LF)-space with
the topology of the inductive limit.

\begin{definition} Let $\Lambda$ be a discrete sequence in $\C$ and let
$\{m_\lambda\}_{\lil}$ be a sequence of natural numbers. The pair
$X=\{(\lambda,m_\lambda)\}_{\lil}$ is called an \emph{interpolating variety} for
the space $A_p$ if for every sequence of values $\{v_\lambda^l\}_{\lambda,l}$,
$\lil$, $l=0,\dots,m_\lambda-1$, with  
\begin{equation}\label{values}
\sup_{\lil} \Bigl(\sum_{l=0}^{m_\lambda-1}|v_\lambda^l|\Bigr) 
e^{-\alpha  p(\lambda)} <\infty 
\end{equation}
for some $\alpha>0$, there exists $f\in A_p$ with 
\[
\frac{f^{(l)}(\lambda)}{l!}=v_\lambda^l\qquad \lil\ ;\quad l=0,
\dots,m_\lambda-1\ .
\]
\end{definition}
The choice of condition \eqref{values} on the values to interpolate reflects
the fact that for every $f\in A_p$ there exists $\alpha>0$ such that
\begin{equation*}
\sup_{z\in\C} \Bigl(\sum_{l=0}^{\infty}\Bigl| 
\frac{f^{(l)}(z)}{l!}\Bigr|\Bigr) 
e^{-\alpha  p(z)} <\infty .
\end{equation*}
Thus, denoting by $A_p(X)$ the space of sequences 
$\{v_\lambda^l\}_{\lambda,l}$ satisfying \eqref{values} for some $\alpha>0$, we
can equivalently define interpolating varieties $X$ as those such that the
restriction operator
\[ 
\begin{split}
\mathcal R_X : A_p &(\C)\longrightarrow A_p(X)\\
& f\quad \mapsto\; \Bigl\{\frac{f^{(l)}(\lambda)}{l!}\Bigr\}_{\lambda,l}
\end{split} 
\]
is onto.

There exists an analytic characterization of $A_p$-interpolating varieties for
general subharmonic  weights $p$ (see \cite[Corollary 3.5]{Br-Li}).  Results
for $p(z)=|z|$ and the weight \eqref{pes} were previously obtained respectively
by Leont'ev \cite{Le} and Squires \cite{Sq0}.  

Given a holomorphic function $f$ let $\mathcal Z(f)$ denote its zero variety,
i.e, the set of pairs $(z,m_z)\in\C\times\N$ such that $f(z)=0$ with
multiplicity $m_z$.

\begin{theoremA}\label{analytic} A variety $X=\{(\lambda,m_\lambda)\}_{\lil}$ is
$A_p$-interpolating if and only if there exists $f\in A_p$ such that
$X\subset \mathcal Z(f)$ and for some constants
$\delta,C>0$
\[
\Bigl|\frac{f^{(m_\lambda)}(\lambda)}{m_\lambda!}\Bigr|
\geq \delta e^{-C  p(\lambda)}\qquad \lil.
\]
\end{theoremA}

We would like to give a geometric description of $A_p$-interpolating varieties
for the non-iso\-tropic Beurling weights
\[
p(z)=|\Im\, z|+\omega(|z|),
\]
where $\omega(t)$ is a subadditive increasing continuous function, 
normalized with $\omega(0)=0$ and such that:
\begin{itemize}
\item[(W1)] $\log(1+t)\lesssim\omega(t)$ for $t>1$.
\item[(W2)] $\displaystyle\int_0^\infty \frac{\omega(t)}{1+t^2} dt<\infty$.
\end{itemize}
Canonical examples of such weights are given by $\omega(t)=\log(1+t^2)$ and
$\omega(t)=t^\gamma$, $\gamma\in(0,1)$.

Beurling weights satisfy the following additional properties:
\begin{itemize}
\item[(c)] For every $c>0$ there exists  $C>0$ such that  $p(\zeta)\leq Cp(z)\;
$ if $\zeta\in D(z,c p(z))$.
\item[(d)] For $\varepsilon>0$ small enough, 
there exists $C(\varepsilon)>0$ such that
if $z\in D(\zeta,\varepsilon p(\zeta))$ then $p(\zeta)\leq C(\varepsilon)
p(z)$. Also, $C(\varepsilon)$ tends to 1  as $\varepsilon$ goes to 0.
\item[(e)] For $x\in\R^+$ big enough, the function $\omega(x)$ does not oscillate too
much. More precisely, fixed $C>0$, if $y\in (x-C\omega(x),x+C\omega(x))$ then 
$1/2\le \omega(y)/\omega(x)\le 2$ for $x$ big enough.  
\end{itemize}
Properties (c) and (d) follow easily from the subadditivity of $\omega$.
Property (e) follows from the subadditivity and the fact that
$\omega(x)=o(|x|/\log|x|)$ (see \cite[Lemma~1.2.8]{Bjorck66}): for any 
$y\in(x-C\omega(x),x+C\omega(x))$
\[
\begin{split}
\omega(x-C\omega(x))&\le \omega(y)\le \omega(x+C\omega(x))\le 
\omega(x-C\omega(x))+\omega(2C \omega(x))\le
 \\ &\le \omega(x-C\omega(x)) + \omega(2Cx/\log x)\le 2 \omega(x-C\omega(x)).
\end{split}
\]
In order to state the geometric conditions on a variety $X$ as above, we
consider the counting function $n(z,r)=\sum\limits_{\lambda\in D(z,r)}
m_\lambda$ and the integrated version
\[
N(z,r)=\int_0^r\frac{n(z,t)-n(z,0)}t\, dt + n(z,0)\log r .
\]
In case we want to specify the variety $X$ which the functions $n$ and  $N$ are
referred to, we will use the notation $n(z,r,X)$ and $N(z,r,X)$ respectively.

We are ready to state our main result. 

\begin{theorem}\label{teorema}
A variety $X=\{(\lambda, m_\lambda)\}_{\lil}$ is $A_p$-interpolating if 
and only if:
\begin{itemize}
\item[(a)] There is $C>0$ such that 
\begin{equation*}
N(\lambda,  p(\lambda), X ) \le C p(\lambda)\qquad \forall 
\lambda\in\Lambda.
\end{equation*}

\item[(b)] The following Carleson-type condition holds
\[
\sup_{x\in\R}
\sum\limits_{\lambda: |\Im\, \lambda|>\omega(|\lambda|)}
m_{\lambda}
\frac{|\Im\, \lambda|}{|x-\lambda|^2}
<\infty .
\]
\end{itemize}
\end{theorem}

Since the Poisson kernel at $\lambda$ in the corresponding half-plane 
(upper half-plane if $\Im\,\lambda>0$ and lower half-plane when 
$\Im\,\lambda<0$) is $P(\lambda,x)=\frac{|\Im\, \lambda|}{|x-\lambda|^2}$, 
a restatement of condition (b) is that the measure 
$\sum\limits_{\lambda: |\Im\, \lambda|>\omega(|\lambda|)}m_\lambda\delta_\lambda$
has bounded Poisson balayage.

\begin{remark} \label{h} 
Notice that for sequences $\Lambda$ within the region
$|\Im\, z|\leq \omega(|z|)$ , condition (a) (shown to be
necessary by Ehrenpreis \& Malliavin and Squires) provides a complete
characterization. However, this is not the case in general, i.e. condition (b)
does not follow from (a), as it is shown in the following example. Take the
sequence $\Lambda$ contained in the angle $\mathcal A=\{z\in\C;\ |\Re z|<\Im
z\}$, having in each segment $\{\Im z=2^n\}\cap \mathcal A$ exactly $2^n$
equispaced points.  Then $\Lambda$ satisfies condition (a) (basically
$n(\lambda,t)\le t$ for $t\le p(\lambda)$), but it does not satisfy (b) (it is
not even a Blaschke sequence).
\end{remark}

\section{Necessary Conditions} 
A standard feature of the $A_p$ spaces is that the interpolation  can be
performed in a stable way. This is a consequence of the open mapping theorem
for (LF)-spaces  applied to the restriction mapping $\mathcal R_X$ defined
in~\ref{preliminaries} (see \cite[Lem\-ma 2.2.6]{Br-Ga}).

\begin{lemma}\label{peak}
If $X$ is an interpolating variety,
there exist $C>0$, $M\in\N$ such that  for every $\lil$ there are functions
$f_\lambda, g_\lambda\in A_p$ with  bounded norms 
$\|f_\lambda\|_{A_{p,M}^\infty},\
\|g_\lambda\|_{A_{p,M}^\infty}\leq C$ and 
\[ 
\begin{split} 
f_\lambda^{(l)}(\lambda')/l!&=\delta_{\lambda\lambda'} 
\delta_{l0}\\
g_\lambda^{(l)}(\lambda')/l!&=\delta_{\lambda\lambda'} \delta_{l(m_\lambda-1)}
\qquad \forall \lambda,\lambda'\in \Lambda,\
0\le l\le m_\lambda.
\end{split} 
\]
\end{lemma}

An application of Jensen's formula to the functions $f_\lambda$, $g_\lambda$ in
the disk $D(\lambda, p(\lambda))$ gives the following result (see \cite[Theorem
4]{Eh-Ma} or \cite[Theorem 1]{Sq}).

\begin{theorem} 
If $X$ is $A_p$-interpolating then condition (a) of Theorem~\ref{teorema}
holds.
\end{theorem}

The necessity of condition (b) 
is an immediate consequence of the following result. 
Assume that $\Lambda\cap \R=\emptyset$; otherwise move the horizontal line so
that it does not touch any of the points in $\Lambda$. Let $\H$ denote the
upper half-plane.

\begin{proposition}\label{BU}
Let $X$ be $A_p$-interpolating. There exist $C>0$ such that
\[
\sum_{
\begin{subarray}{c}
\lambda'\in\Lambda\cap\H\\
\lambda'\neq\lambda 
\end{subarray}
}
m_{\lambda'} \log\Bigl|\frac{\lambda-\lambda'}{\lambda-\bar\lambda'}\Bigr|^{-1}
\leq C p(\lambda)\quad \textrm{for all $\lil\cap\H$.} 
\]
\end{proposition}

Of course an analogous result could be given for any upper ($\{z: \Im z>a\}$) 
or lower ($\{z: \Im z<a\}$) half plane.

In order to see that Proposition~\ref{BU}  implies condition (b) of
Theorem~\ref{teorema} define $\Lambda_+=\Lambda\cap \{\Im\, z>\omega(|z|)\}$.
Given $x\in\R$ consider $\lambda\in\Lambda_+$ such that 
$|x-\lambda|=\inf_{\Lambda_+} |x-\lambda|$. Then
\[
|\lambda-\bar\lambda'|\leq |\lambda-x| + |x-\bar\lambda'|=
|\lambda-x| + |x-\lambda'|\leq 2 |x-\lambda'|,
\]
and therefore
\[
\sum_{\lambda'\in\Lambda_+} m_{\lambda'} 
\frac{|\Im \, \lambda'|}{|x-\bar\lambda'|^2}
\leq 2\sum_{\lambda'\in\Lambda_+} m_{\lambda'} 
\frac{|\Im \, \lambda'|}{|\lambda-\bar\lambda'|^2}.
\]
The estimate $\log t^{-1}\geq 1-t$ for $t\in(0,1)$ shows that
\begin{equation*}
\sum_{
\begin{subarray}{c}
\lambda'\in\Lambda_+\\
\lambda'\neq\lambda 
\end{subarray}
}
m_{\lambda'} \frac{|\Im \, \lambda| |\Im \, \lambda'|}{|\lambda-\bar\lambda'|^2}
\leq
\sum_{
\begin{subarray}{c}
\lambda'\in\Lambda_+\\
\lambda'\neq\lambda 
\end{subarray}
}
m_{\lambda'} \log\Bigl|\frac{\lambda-\lambda'}{\lambda-\bar\lambda'}\Bigr|^{-1}.
\end{equation*}

Since $p(\lambda)\simeq |\Im \, \lambda|$ for $\lambda\in\Lambda_+$, it is
clear that this implies condition (b) of Theorem~\ref{teorema}.

\begin{remark}\label{7}
The necessary condition of Proposition~\ref{BU} can be seen as a Carleson
type condition; it can be rewritten as
\[ 
\begin{split}
|B_\lambda(\lambda)|&\geq \delta e^{-Cp(\lambda)}\quad \lambda\in\Lambda\cap\H,
\end{split} 
\]
where $B$ denotes the Blaschke product in $\H$ of
$\{(\lambda,m_\lambda)\}_{\lambda\in\Lambda\cap\H}$, and 
\[
B_\lambda(z)=
B(z)\Bigl(\frac{z-\bar\lambda} {z-\lambda}\Bigr)^{m_\lambda}.
\]

It can also be seen as density conditions for the counting function associated
to the hyperbolic metric in the half-plane. Letting
$\nu=\sum_{\lambda\in\Lambda\cap\H} m_\lambda \delta_\lambda$  and using the
distribution function we have
\[ 
\begin{split}
\sum_{\lambda\in\Lambda\cap\H}
m_\lambda
\log\Bigl|\frac{z-\lambda}{z-\bar\lambda}\Bigr|^{-1}&=\int_{\H}
m_\lambda\log\Bigl|\frac{z-\zeta}{z-\bar\zeta}\Bigr|^{-1} d\nu(\zeta)
=\int_0^1 \frac {n_{\H}(z,t)}t\;  dt  ,
\end{split} 
\]
where
\[
D_{\H}(z,t)=
\{\zeta : \Bigl|\frac{z-\zeta}{z-\bar\zeta}\Bigr|<t\} 
\quad\textrm{and}\quad
n_{\H}(z,t):=\nu(D_{\H}(z,t))
\]
is the number of points of $\Lambda$ in the pseudohyperbolic disk of ``center"
$z$ and ``radius" $t$ (actually the true disk of center $\Re z +i
\frac{1+t^2}{1-t^2}\Im z$ and radius $\frac{2t}{1-t^2}\Im z$).
\end{remark}

\begin{proof}
Let $z=x+iy$ and consider the Poisson transform of $\omega(|t|)$:
\[
u(z):=P[\omega](z)=\int_{\R}\frac{y\, \omega(|t|)}{(x-t)^2+y^2}\, dt,
\]
which converges by (W2). Define $H=\exp(u+i\tilde u)$, where $\tilde u$
is a harmonic conjugate of $u$.

Given $\lambda\in\Lambda\cap \H$, take the function $f_\lambda$ given by
Lemma~\ref{peak}  and define
\[
h_\lambda(z)=\frac{f_\lambda(z) e^{iM_1z}}{H^{M_2}(z)} ,
\]
with $M_1,M_2$ to be chosen.
It is clear that $h_\lambda$ is holomorphic in $\H$. On the other hand, for all $z$ in the upper half plane
$|\log|H(z)|- \omega(|\Re z|)|\le A + B|\Im z|$, 
see \cite[Lemma~1.3.11]{Bjorck66}. Moreover $|\omega(|\Re z|)-\omega(|z|)|\le
\omega(|\Im z|)\le A+B|\Im z|$, thus 
$|\log|H(z)|- \omega(|z|)|\le A + B|\Im z|$.
Therefore, if $M_1$ and $M_2$ are big enough,
$h_\lambda$ is bounded in $\H$ by a constant which does 
not depend on $\lambda$:
\[
|h_\lambda(z)|\leq C e^{M p(z)-M_1\Im\, z-M_2\log|H(z)|}\lesssim 1 .
\]
Also, 
\[
|h_\lambda(\lambda)|=e^{-M_1\Im\, \lambda-M_2\log|H(\lambda)|}
\geq e^{-C p(\lambda)}.
\]
Apply now Jensen's formula in the half-plane to the function $h_\lambda$:
\[
\log|h_\lambda(\lambda)|=\int_{\R} P(\lambda,x) \log |h_\lambda(x)|\, dx
-\int_{\H} G(\lambda,\zeta)\Delta \log |h_\lambda(\zeta)|,
\]
where $P(\lambda,x)$ denotes the Poisson kernel and
$G(\lambda,\zeta)=
\displaystyle\log\Bigl|\frac{\lambda-\zeta}{\lambda-\bar\zeta}\Bigr|^{-1}$
is the Green function in $\H$ with pole in $\lambda$.

Since $h_\lambda$ vanishes on $\Lambda\setminus\{\lambda\}$, this and  
the estimates above yield
\[
\sum_{
\begin{subarray}{c}
\lambda'\in\Lambda\cap\H\\
\lambda'\neq\lambda 
\end{subarray}
}
m_{\lambda'}
\log\Bigl|\frac{\lambda-\lambda'}{\lambda-\bar\lambda'}\Bigr|^{-1}
\leq \sup_{\R} \log|h_\lambda|-\log|h_\lambda(\lambda)|\lesssim p(\lambda).
\]
\end{proof}

\section{Sufficient conditions}
We split the sequence into three pieces, according to the non-isotropy of the
weight $p$. Consider the regions
\[ 
\begin{split}
\Omega_0&=\{z\in\C : |\Im\, z|\leq \omega(|z|) \}\\
\Omega_+&=\{z\in\C : \Im\, z > \omega(|z|)\}\\
\Omega_-&=\{z\in\C : \Im\, z < -\omega(|z|)\} ,
\end{split} 
\]
and define $\Lambda_0=\Lambda\cap\Omega_0$,  $\Lambda_+=\Lambda\cap\Omega_+$
and $\Lambda_-=\Lambda\cap\Omega_-$. Let also
$X_0=\{(\lambda,m_\lambda)\}_{\lambda\in\Lambda_0}$, 
$X_+=\{(\lambda,m_\lambda)\}_{\lambda\in\Lambda_+}$ and
$X_-=\{(\lambda,m_\lambda)\}_{\lambda\in\Lambda_-}$.

It is enough to prove that each piece $X_+,X_-,X_0$ of the variety $X$ is
$A_p$-in\-ter\-po\-la\-ting. This is so because $X$ is weakly
separated (see Lemma~\ref{lemma} (a) below), and a weakly separated union of a
finite number of $A_p$-interpolating varieties is also  $A_p$-interpolating
\cite[Theorem II.1]{Ou03}. It is also clear that the varieties $X^+$ and
$X^-$ can be dealt similarly.

We start with some easy consequences of condition (a) of Theorem~\ref{teorema}.

\begin{lemma}\label{lemma}
If condition (a) in Theorem~\ref{teorema} holds, then 
\begin{itemize}
\item[(a)] $X$ is weakly separated: there exist
$\delta,C>0$ such that the disks $D_\lambda=D(\lambda, \delta
e^{-C\frac{p(\lambda)}{m_\lambda}})$ are pairwise disjoint, i.e.
\[
|\lambda-\lambda'|\geq 2\delta \max[ e^{-C\frac{p(\lambda)}{m_\lambda}}, 
e^{-C\frac{p(\lambda')}{m_{\lambda'}}}]
\qquad \lambda\neq \lambda'.
\]

\item[(b)] There exist  $\varepsilon,C>0$
such that $n(z,\varepsilon p(z), X)\le C p(z)$, $\forall z\in \C$.
\end{itemize}
\end{lemma}

\begin{proof}
(a)  If there exists $\lambda'$ such that $|\lambda'-\lambda|<1$ then
\[
N(\lambda,p(\lambda),X)\geq\int_{|\lambda'-\lambda|}^1
\frac{n(\lambda,t)-m_\lambda}t dt\geq
\int_{|\lambda'-\lambda|}^1\frac{m_{\lambda'}}t dt=
\log\Bigl(\frac 1{|\lambda'-\lambda|}\Bigr)^{m_{\lambda'}}.
\]
Using Theorem~\ref{teorema} (a) and reversing the roles of $\lambda$ and
$\lambda'$ we obtain the desired estimate. 

(b) When $z=\lambda\in\Lambda$, this is
immediate from the estimate
\[
\int_{ 1/2 p(\lambda)}^{p(\lambda)} \frac{n(\lambda, 
1/2 p(\lambda))-1}t dt\leq  N(\lambda,p(\lambda)) .
\]
When $z\notin \Lambda$, then  let $\varepsilon>0$ be such that $\zeta\in
D(z,\varepsilon p(z))$ implies  
\[ 
D(z, \varepsilon p(z))\subset  D(\zeta,1/2 p(\zeta)),
\]
which exists by property (d) of the weight. Take $\lambda\in
D(z,\varepsilon p(z))$ (if there is no such $\lambda$ the estimate is
obviously  true). Then, by the previous case and property (c) of the weight
\[
n(z,\varepsilon  p(z))\leq n(\lambda, 1/2 p(\lambda))\lesssim p(\lambda)
\lesssim p(z).
\]
\end{proof}

\subsection{ Case $\Lambda_0$}  We would like to prove that
$X_0=\{(\lambda,m_\lambda)\}_{\lambda\in\Lambda_0}$ is $A_p$-interpolating
using a $\bar\partial$-scheme.  This is easier if we can regularize the weight
in the following way.

\begin{lemma}\label{regular}
There exists $\tilde p$ subharmonic in $\C$ such that $p(z)\simeq \tilde p(z)$
and
\begin{equation}\label{regularization}
%1/p(z)\simeq 
1/ \tilde p(z)\lesssim \Delta\tilde p(z) \qquad\textrm{if} \quad 
|\Im z|\leq 2\omega(|z|).
\end{equation}
\end{lemma}

The fact that $p\simeq \tilde p$ clearly implies that 
$A_p=A_{\tilde p}$ and the interpolating varieties for
$A_p$ and $A_{\tilde p}$ are the same.

\begin{proof}
We will 
construct $\tilde p(z) = |\Im z|+r(z)$, where $r$ satisfies the following
properties:
\begin{itemize}
\item[(i)] $r\geq 0$  and $\tilde p$ is subharmonic in $\C$,
\item[(ii)]$r(z)=0$ if $|\Im z|\geq 10\omega(|z|)$.
\item[(iii)]   $1/p(z)\lesssim \Delta \tilde p(z)$ 
and $r(z)\simeq \omega(|z|)$ if $|\Im z|\leq 2\omega(|z|)$.
\end{itemize}
In order to construct $r$ partition the real line in intervals $I_n$ of
center $x_n$ and length  $\omega_n=\omega(x_n)$. 

We consider two measures in $\C$. The first one is the usual length measure
$d\nu$ in $\R$, which we split $d\nu=\sum_n d\nu_n$, with $d\nu_n=dx_{|I_n}$. 
The second one is defined as a sum of convolutions of the $d\nu_n$'s: let
\[
d\mu_n(z)=\bigl(\frac1{100\pi\omega_n^2}\int_{I_n} \chi_{D_n}(z-x) dx\bigr)\, dm(z) ,
\] 
where  $D_n=D(0,10\omega_n)$, and define 
$d\mu=\sum_n d\mu_n$.

Notice that when $z$ is at a distance of $I_n$ smaller than $2\omega_n$, we can
use property (e) of the Beurling weights to deduce that $d\mu(z)\simeq
1/\omega(|z|)\simeq1/p(z)$. Hence $d\mu(z)\simeq dm(z)/ p(z)$.

Define 
\[
r(z) = \int_{\C} \log|z-w| (d\mu(w)-d\nu(w)). 
\]
Since $\Delta|\Im\, z|=d\nu$ we have 
$\Delta \tilde p = d\mu\ge 0$. 

Let $S_n$ denote the support of $\mu_n$. Let
\[
r_n(z):=\int_{\C} \log|z-w| (d\mu_n(w)-d\nu_n(w))= 
\int_{S_n}\log|z-w|d\mu_n(w)-\int_{I_n}\log|z-x|dx
\]
Using the definition of $\mu_n$ and reversing the order of 
integration we get
\[
r_n(z)=\int_{I_n}M(x)dx ,
\]
where 
\[
M(x)=\frac 1{100\pi\omega_n^2}\int_{D(x,10\omega_n)}
\log|z-w|dm(w)-\log|z-x|\ge 0.
\]
In particular, $r$ in non-negative in $\C$.

If $z\notin S_n$ and $x\in I_n$, $\log|z-w|$ is harmonic in
$D(x,10\omega_n)$, hence $r_n(z)=0$

Suppose now $z\in D(x_n,3\omega_n)$. Then, for each $x\in I_n$, $|z-x|\leq
4\omega_n$ and
\[
M(x)\geq \frac 1{100\pi\omega_n^2}\int_{9\omega_n\leq|w-x|\leq 10\omega_n}
\log\frac{|z-w|}{|z-x|}dm(w)\gtrsim 1.
\]
Thus, $r_n(z)\gtrsim \omega_n\gtrsim \omega(|z|)$.

If $z\in S_n$, using that $\mu_n$ and $\nu_n$ have the same mass $\omega(x_n)$,
we obtain
\[
\begin{split}
&\int_{\C} \log|z-w| (d\mu_n(w)-d\nu_n(w))\le
\int_{\C} \Bigl|\log\frac{|z-w|}{\omega(x_n)} \Bigr|(d\mu_n(w)+d\nu_n(w))
\lesssim\\
&\int_{\C} \Bigl|\log\frac{|x_n-w|}{\omega(x_n)} \Bigr|(d\mu_n(w)+d\nu_n(w))
\lesssim \omega(|z|).
\end{split}
\]
Since $|\Im z|\leq 2\omega(|z|)$, $z$ belongs at most to a finite  number of
$S_n$'s and at least to  one $D(x_n,10\omega_n)$, by property (e) of the
Beurling weights,  we are done.
\end{proof}

Let us prove now that $X_0$ is $A_{\tilde p}$-interpolating. Denote $p$ instead
of $\tilde p$ and assume that $\Delta p\simeq 1/p$ on  $|\Im z|\leq 2
\omega(|z|)$.

Consider the separation radius $\delta_\lambda:=\delta e^{-Cp(\lambda)}$
given by Lemma~\ref{lemma} (a).

Given a sequence of values $\{v_\lambda\}_{\lil}$ satisfying \eqref{values},
define the smooth interpolating function
\[
F(z)=\sum_{\lil_0} p_\lambda(z) \mathcal X\Bigl(\frac{|z-\lambda|^2}
{\delta^2_\lambda}\Bigr),
\]
where $p_\lambda(z)=\sum\limits_{l=0}^{m_\lambda-1} v_\lambda^l (z-\lambda)^l$
and $\mathcal X$ is a smooth cut-off function with $|\mathcal X^\prime|\lesssim
1$, $\mathcal X(x)=1$ if $|x|\leq 1$ and $\mathcal X(x)=0$ if $|x|\geq 2$.

It is clear that $F^{(l)}(\lambda)/l!=v_\lambda^l$,  and that $F$ has the
characteristic growth of $A_p$ functions: the support of $F$ is contained in
$\cup_\lambda  D_\lambda$ and  for $z\in D_\lambda$
\[
|F(z)|\leq \sum_{l=0}^{m_\lambda-1} |v_\lambda^l| \leq C e^{\alpha p(\lambda)}
\lesssim e^{Kp(z)} .
\]
There is also a good estimate on $\bar\partial F$.
Its support is the reunion of the annuli
\[
C_\lambda=\{z\in\C : \delta_\lambda \leq |z-\lambda|\leq 
2 \delta_\lambda\} ,
\] 
and for $z\in C_\lambda$,
\[
\Bigl|\frac{\partial F}{\partial\bar z}(z)\Bigr|\lesssim 
\sum_{l=0}^{m_\lambda-1} |v_\lambda^l|
|\mathcal X'|\frac1{\delta_\lambda}\lesssim e^{C p(\lambda)}\lesssim e^{Kp(z)},
\]
for $K$ big enough.

Altogether, there exists $\gamma>0$ such that
\begin{equation}\label{estimates}
\int_{\C} |F(z)|^2 e^{-\gamma p(z)}<\infty\quad ,
\quad\int_{\C} |\bar\partial F(z)|^2 e^{-\gamma
p(z)}<\infty\ .
\end{equation}
Now, when looking for a holomorphic interpolating function of the form
$f=F-u$, we are led to the $\bar\partial$-problem
\[
\bar\partial u=\bar\partial F\ , 
\]
which we solve using H\"ormander's theorem \cite[Theorem 4.2.1]{Ho}:
given a (pluri)\-sub\-har\-monic function $\psi$ in $\C$, there exists a
solution $u$
to the above equation such that
\[
2\int_{\C} |u|^2 \frac{e^{-\psi}}{(1+|z|^2)^2}\; dm\leq
\int_{\C} |\bar\partial  F|^2 e^{-\psi} dm \ .
\]
We apply H\"ormander's theorem with
\[
\psi_\beta(z)=\beta  p(z) + v(z)\ ,
\]
where $\beta>0$ will be chosen later on and
\[
v(z)=\sum_{\lil_0} m_\lambda 
\Bigl[\log|z-\lambda|^2-\frac 1{\pi \varepsilon^2 p^2(\lambda)}
\int_{D(\lambda,\varepsilon p(\lambda))}\log|z-\zeta|^2 dm(\zeta)\Bigr]\ .
\]
Here $\varepsilon$ is a fixed small constant to be determined later on.

Integrating by parts the equality 
\begin{equation*}
\int_0^{2\pi}\log|a-re^{i\theta}|^2\; \frac{d\theta}{2\pi}=
\begin{cases}
\log |a|^2 \quad&\textrm{if}\quad |a|> r\\
\log r^2 \quad&\textrm{if}\quad |a|\leq r
\end{cases}
\end{equation*}
one sees that for 
$a\in\C$ and $r>0$:
\begin{equation*}
\log|a|^2-\frac 1{\pi r^2}\int\limits_{D(a,r)}\log|\zeta|^2 dm(\zeta)=
\begin{cases}
\log |\frac ar|^2+1-|\frac ar|^2 \ &\textrm {if}\quad |a|\leq r\\
0 \ &\textrm{if}\quad |a|< r\ .
\end{cases}
\end{equation*}
Thus 
\[
v(z)=\sum_{\lambda:|\lambda-z|\leq\varepsilon p(\lambda)} 
m_\lambda\Bigl[\log\frac{|z-\lambda|^2}{\varepsilon^2 p^2(\lambda)}
+1-\frac{ |z-\lambda|^2}{\varepsilon^2 p^2(\lambda)}\Bigr]\ .
\]
In particular $v\leq 0$ and $\Delta v(z)=0$ if 
$z\notin\cup_\lambda D(\lambda,\varepsilon p(\lambda))$.
For $z\in\cup_\lambda D(\lambda,\varepsilon p(\lambda))$ we have
$|\Im\, z|\leq 2\omega(|z|)$ and 
\[
\Delta v(z)\geq \! \sum_{\lambda:|\lambda-z|\leq\varepsilon  p(\lambda)}   
\frac {-m_\lambda}{\varepsilon^2 p^2(\lambda)}\gtrsim \!
\sum_{\lambda:|\lambda-z|\leq C(\varepsilon) p(z)} 
\frac {-m_\lambda}{p^2(z)}
= - \; \frac{n(z, C(\varepsilon) p(z))}{p^2(z)}\ .
\]
As observed in Lemma~\ref{lemma} (b), with $\varepsilon$ small enough 
$n(z,C(\varepsilon) p(z))\lesssim p(z)$, thus $\Delta v(z)\gtrsim -1/p(z)$.
This and \eqref{regularization} show that $\psi_\beta$ is subharmonic if
$\beta$ is chosen big enough.

Also, we deduce from (W1) that for any $\beta'>\beta$:
\[
\int_{\C} |u|^2 e^{-\beta'p }dm\lesssim \int_{\C} |u|^2 
\frac{e^{-\psi_\beta}}{(1+|z|^2)^2}\; dm\lesssim
\int_{\C} |\bar\partial F|^2 e^{-\psi_\beta} dm\ .
\]

We need to control $\psi_\beta$ on the support of $\bar\partial F$. 
For $z\in C_\lambda$,
\[ 
\begin{split}
|\psi_\beta(z)-\beta p(z)|&\leq 
\sum_{\lambda :|\lambda-z|\leq \varepsilon p(\lambda)}
m_\lambda  
\log \frac {\varepsilon^2 p^2(\lambda)}{|z-\lambda|^2}\\
&\simeq m_\lambda 
\log \frac {\varepsilon^2 p^2(\lambda)}{|z-\lambda|^2} +\sum_{
\begin{subarray}{c}
\lambda': |z-\lambda'|\leq\varepsilon p(\lambda')\\
\lambda'\neq \lambda
\end{subarray}
}
m_{\lambda'}
\log \frac {\varepsilon^2 p^2(\lambda')}{|z-\lambda'|^2}\\
&\lesssim p(\lambda) + \sum_{
\begin{subarray}{c}
\lambda': |\lambda'-z|\leq C(\varepsilon) p(z)\\
\lambda'\neq \lambda
\end{subarray}
}
m_{\lambda'}
\log \frac { C(\varepsilon)^2 p^2(z)}{|z-\lambda'|^2}
\lesssim p(z) + N(z,C(\varepsilon) p(z))
\end{split} 
\]

\begin{claim}\label{claim}
For $\varepsilon$ small enough $N(z,C(\varepsilon) p(z))\lesssim
p(z)$ for all $z\in supp(\bar\partial F)$. 
\end{claim}
Assuming the claim we have 
$|\psi_\beta(z)-\beta p(z)| \leq K p(z)$ on $supp(\bar\partial F)$.
Therefore, for $\beta$ big enough 
\[
\int_{\C} |u|^2 e^{-\beta'p }dm \lesssim \int_{\C}
|\bar\partial F|^2 e^{-\psi_\beta} dm\leq \int_{\C} |\bar\partial 
F|^2 e^{-\gamma p }
dm<\infty\ . 
\]
This shows that $f:=F-u\in A_p$. Since $e^{-\psi_\beta}\simeq
|z-\lambda|^{-2m_{\lambda}}$   around each $\lambda$, also $u^{(l)}(\lambda)=0$
for all $\lil$, $l=0,\dots,m_{\lambda}-1$,  and therefore $
f^{(l)}(\lambda)/l!=F^{(l)}(\lambda)/l!=v_\lambda^l$, as required. 

\textsl{Proof of the claim:}  Assume $z\in C_\lambda$ and observe that
$n(z,t)=0$ for $t<\delta_\lambda$ and that $n(z,t)\leq m_\lambda$ for 
$\delta_\lambda\leq t< 2\delta_\lambda$. Since $D(z,t)\subset D(\lambda,
t+2\delta_\lambda)$ and  $|z|<|\lambda|+2 \delta_\lambda$, we have (changing
into $s=t+2\delta_\lambda$)
\[ 
\begin{split}
N(z,C(\varepsilon) p(z))&\leq\int_{\delta_\lambda}^{2\delta_\lambda} 
\frac {m_{\lambda}}t\; dt +
\int_{2\delta_\lambda}^{C(\varepsilon) p(z)}
\frac{n(z,t)-m_\lambda}t \; dt \leq \\
&\leq p(\lambda) + \int_{4\delta_\lambda}^{C(\varepsilon) p(z)+2\delta_\lambda}
\frac{n(\lambda,s)-m_\lambda}{s-2\delta_\lambda} \; ds\lesssim\\
&\lesssim  p(\lambda) + \int_{4 \delta_\lambda}^{C(\varepsilon) p(z) 
+2\delta_\lambda}
\frac{n(\lambda,s)-m_\lambda}{s/2} \; ds
\lesssim  p(\lambda)+ 
N(\lambda,C'(\varepsilon) p(\lambda))\ .  
\end{split} 
\]

From the properties of the weight and the hypothesis we have finally that for
$\varepsilon$ small\newline 
$N(z,C(\varepsilon) p(z))\lesssim
p(\lambda) \lesssim p(z)$. 

\subsection{Case $\Lambda^+$}  According to Theorem A, it
is enough to construct a function  $G\in A_p$  such that $X_+\subset \mathcal
Z(G)$ and
\[
\frac{|G^{(m_\lambda)}(\lambda)|}{m_\lambda!}\geq \varepsilon e^{-K p(\lambda)}
\qquad \lambda\in\Lambda_+ 
\]
for some constants $\varepsilon, k>0$. In fact, the hypotheses of Theorem~A
require the weight $p$ to be subharmonic, and our weights are not necessarily
so. Nevertheless, by Lemma~\ref{regular}, there exists a  subharmonic weight
$\tilde p$ equivalent to $p$, and we may apply Theorem~A to $\tilde p$.

Take any entire function $F$ such that $\mathcal Z(F)=X_+ $. Since the
necessary conditions imply that $X_+$ satisfies the Blaschke condition in
$\H$,  we can consider also the Blaschke product 
\[
B(z)=\prod_{\lambda\in\Lambda_+}\Bigl(\frac{z-\lambda}
{z-\bar\lambda}\Bigr)^{m_\lambda}, \qquad z\in\H .
\]
Define
\begin{equation*}
\phi(z)=
\begin{cases}
\displaystyle\log\Bigl|\frac{F(z)}{B(z)}\Bigr|\quad &\textrm{ $\Im\, z>0$}\\
\log|F(z)|\quad &\textrm{ $\Im\, z\leq0$.}
\end{cases}
\end{equation*}

\begin{lemma}
$\phi\in\Har(\C\setminus\R)\cap  SH(\C)$ and its Laplacian
is uniformly bounded.
\end{lemma}

\begin{proof} It is clear, by definition, that $\phi\in\Har(\C\setminus\R)$.
In order to prove that $\phi\in SH(\C)$ it is enough to check the mean
inequality for  $x\in\R$. We have
\[
\phi(x)=\log |F(x)|\leq\frac 1{2\pi}\int_0^{2\pi} \log |F(x+r e^{i\theta})| 
d\theta\leq\frac 1{2\pi}\int_0^{2\pi} \phi(x+r e^{i\theta}) d\theta .
\]
Since
$\Delta\log|F|\equiv 0$ around $\R$, it is enough to compute the
Laplacian of
\begin{equation*}
\psi(z)=
\begin{cases}
\displaystyle\log\frac{1}{|B(z)|}\quad &\textrm{ $\Im\, z>0$}\\
\quad 0\quad &\textrm{ $\Im\, z\leq0$}.
\end{cases}
\end{equation*}
Being
\[
\log \frac{1}{|B(z)|}=\frac 12
\sum_{\lambda\in\Lambda^+} m_\lambda 
\log \Bigl|\frac{z-\bar\lambda}{z-\lambda}\Bigr|^2 ,
\]
it will be enough to compute the Laplacian of each term
\begin{equation*}
\psi_\lambda(z)=
\begin{cases}
\displaystyle\log \Bigl|\frac{z-\bar\lambda}{z-\lambda}\Bigr|^2\quad 
&\textrm{ $\Im\, z>0$}\\
\quad 0\quad &\textrm{ $\Im\, z\leq0$}.
\end{cases}
\end{equation*}

It is clear that $\partial \psi_\lambda/\partial x=0$ on $\R$, hence $\Delta
\psi_\lambda= \partial^2 \psi_\lambda/\partial y^2$. Since $\psi_\lambda$ is
continuous around $\R$, this Laplacian has a magnitude equivalent to the jump
of the first derivative of $\psi_\lambda$. The derivative of the Green 
function on the
half-plane with respect to the normal direction $y$ is the Poisson kernel:
\[
\frac{\partial}{\partial y}\log\Bigl|\frac{z-\bar\lambda}{z-\lambda}
\Bigr|^2_{|y=0}=
\frac{4\Im\, \lambda}{|x-\lambda|^2}.
\]
Therefore
\[
\Delta\phi(x)=4
\displaystyle\sum_{\lambda\in\Lambda^+} m_\lambda 
\frac{\Im\,\lambda}{|x-\lambda|^2}\,dx,
\]
which is bounded by hypothesis.
\end{proof}

Define
\[
\Psi(z)=N|\Im\, z|-\phi(z).
\]

Observe that $\Delta\Psi(z)=N\, dx - \Delta\phi(x)\, dx$, thus according to
the previous Lemma $\Delta\Psi\simeq dx$ when  $N\in\N$ is big enough. In this
situation, according to \cite[Lemma~3]{Or-Se},  there exists a multiplier
associated to $\Psi$, i.e., an entire function $h$ such that:
\begin{itemize}
\item[(a)] $\mathcal Z (h)$ is a separated sequence contained in $\R$
\item[(b)] Given any $\varepsilon>0$, 
$|h(z)|\simeq \exp(\Psi(z))$ for all points $z$ such that $d(z,\mathcal Z(h))
>\varepsilon$.
\end{itemize}

Define now $G= h F$. It is clear that $G\in A_p$:
\[
|G(z)|\lesssim e^{\Psi(z)+\log |F(z)|}\leq e^{\Psi(z)+\phi(z)}\leq 
e^{N p(z)}\qquad z\in\C .
\]
It is also clear that $X_+\subset\mathcal Z(G)$, since
$X_+\subset\mathcal Z(F)$. 

In order to prove that there exist $\varepsilon,C>0$ such that
\begin{equation}\label{cota} 
\Bigl|\frac{G^{(m_\lambda)}(\lambda)}{m_\lambda!}\Bigr|\geq\varepsilon
e^{-C p(\lambda)}
\end{equation}
consider then the disjoint disks $D_\lambda=D(\lambda,\delta_\lambda)$,
$\delta_\lambda=\delta e^{-C\frac{p(\lambda)}{m_\lambda}}$
given by Lemma~\ref{lemma}(a).
Since $\Lambda_+$ is far from $\mathcal Z(h)$, the estimate
\[
|G(z)|=|h(z)|e^{\phi(z)}|B(z)|\simeq e^{N|\Im\, z|}  |B(z)| 
\qquad z\in\partial D_\lambda
\]
holds.

\begin{claim}
There exists $C>0$ such that $|B(z)|\geq\epsilon e^{-C p(z)}$, 
$z\in\partial D_\lambda$.
\end{claim}

Assuming this we have $|G(z)| \gtrsim  e^{-C p(z)} $ for all 
$z\in\partial D_\lambda$.
Define then $g(z)=G(z)/(z-\lambda)^{m_\lambda}$. It is clear that $g$ is
holomorphic, non-vanishing in $D_\lambda$, and  $|g(z)|\gtrsim e^{-c
p(\lambda)}$ for $z\in\partial D_\lambda$. By the minimum principle 
\[
\Bigl|\frac{G^{(m_\lambda)}(\lambda)}{m_\lambda !}\Bigr|=|g(0)|\gtrsim e^{-c
p(\lambda)},
\] 
as desired.

\textsl{Proof of the claim:} As observed in Remark~\ref{7}(b), the estimate we
want to prove is equivalent to 
\[
\int_0^1\frac{ n_{\H}(z,t)}t \; dt\lesssim p(z)\qquad 
z\in\partial D_\lambda .
\]
This is proved as Claim~\ref{claim}, replacing the Euclidean disks by
the hyperbolic ones.
We have
\[
\int_0^1\frac{ n_{\H}(z,t)}t \; dt\lesssim 
\int_{\delta_\lambda}^{2\delta_\lambda}\frac{m_\lambda}t dt  +
\int_{2\delta_\lambda}^1\frac{ n_{\H}(z,t)-m_\lambda}t \; dt.
\]
The first term is controlled by $p(\lambda)$. In
order to control the second term observe that $D_{\H}(z,t)\subset
D_{\H}(\lambda,\frac{t+\delta_\lambda}{1+t \delta_\lambda})$; hence changing
the variable into $s=\frac{t+\delta_\lambda}{1+t \delta_\lambda}$ we get
\[ 
\begin{split}
\int_{2\delta_\lambda}^1\frac{ n_{\H}(z,t)-m_\lambda}t \; dt\leq
\int_{\frac{3\delta_\lambda}{1+2\delta_\lambda^2}}^1
\frac{ n_{\H}(\lambda,s)-m_\lambda}{s-\delta_\lambda} 
\frac{1-\delta_\lambda^2}{(1-\delta_\lambda)^2}\; ds .
\end{split} 
\]
There is no restriction in assuming that $\delta_\lambda<1/2$. Then
$\frac{3\delta_\lambda}{1+2\delta_\lambda^2}>2\delta_\lambda$ and therefore
$s-\delta_\lambda>s/2$. With this and Theorem~\ref{teorema}(b) 
we obtain
\[
\int_0^1\frac{ n_{\H}(z,t)}t \; dt\lesssim p(\lambda)+
\int_0^1 \frac{ n_{\H}(\lambda,s)-m_\lambda}s \; ds\ .
\]
Since $p(\lambda)\lesssim p(z)$, we will be done as soon as we prove that
\[
\int_0^1 \frac{ n_{\H}(\lambda,s)-m_\lambda}s \; ds\lesssim p(\lambda).
\]
There exists $\delta>0$ (independent of $\lambda$) such that
$D_{\H}(\lambda,\delta)\subset D(\lambda,p(\lambda))$. Then
\[ 
\begin{split}
\int_0^\delta \frac{ n_{\H}(\lambda,s)-m_\lambda}s \; ds &
=\!\!\!\sum\limits_{0<|\frac{\lambda-\lambda'}{\lambda-\bar\lambda'}|<\delta}
\!\!m_{\lambda'}\log\frac{\delta}{|\frac{\lambda-\lambda'}{\lambda-
\bar\lambda'}|}\leq\!\!\!\!
\sum\limits_{0<|\frac{\lambda-\lambda'}{\lambda-\bar\lambda'}|<\delta }
m_{\lambda'}\log\frac{p(\lambda)}{|\lambda-\lambda'|}\\
&\lesssim\!\!\!\!
\sum\limits_{0<|\lambda-\lambda'|<p(\lambda)}\!\!
m_{\lambda'}\log\frac{p(\lambda) }{|\lambda-\lambda'|}\leq
N(\lambda,p(\lambda))\lesssim p(\lambda).
\end{split} 
\]

For the remaining part we use Theorem~\ref{teorema}(b) 
and the estimate $\log
t^{-1}\simeq 1-t$ for $\delta<t< 1$. Taking $x=\Re\,  \lambda$ we have
\[
\begin{split}
\int_\delta^1 \frac{ n_{\H}(\lambda,s)-m_\lambda}s \; ds&\lesssim
\sum_{\lambda\neq\lambda'} m_{\lambda'}
\frac{|\Im\, \lambda| |\Im\, \lambda'|}{|\lambda-\bar\lambda'|^2}\lesssim\\
&\lesssim
\sum_{\lambda\neq\lambda'} m_{\lambda'}
\frac{|\Im\, \lambda| |\Im\, \lambda'|}{|x-\lambda'|^2}\lesssim
|\Im\, \lambda|\simeq p(\lambda).
\end{split}
\]

\providecommand{\bysame}{\leavevmode\hbox to3em{\hrulefill}\thinspace}
\providecommand{\MR}{\relax\ifhmode\unskip\space\fi MR }

\end{document}